\documentstyle[11pt, amsfonts,epsf,here ]{article}

\def\C{{\mathbb C}}

\def\a{\alpha}
\def\b{\beta}
\def\e{\varepsilon }

\def\n{\newline}
\def\w{\newline \newline}

\begin{document}
  
{\bf ON ZEROS OF SOLUTIONS OF LINEAR DIFFERENTIAL SYSTEMS.}\w
\begin{center}{ 
 by {\bf Alexei Grigoriev}}\end{center}
\begin{center}{ \ \ \ \ \ \ \ \ \ \ \ \ \ \ \ \ \ 
Weizmann Institute of Science, Rehovot\w}\end{center}
The subject discussed in this text is the oscillatory behavior of
solutions of holomorphic 
linear time-dependent differential systems of equations.
In [1] a bound is provided for the number of times a trajectory of a
polynomial vector field crosses a given algebraic hypersurface.
This implies a bound also for linear time dependent vector fields with
polynomial dependence on time. A different approach is taken here,
giving better bound for some special cases.\w
This text is a preliminary report.\w
Besides of having an 
intrinsic interest, a related problem arises in estimating the number of
zeros which the so called complete Abelian integrals may have, the
latter problem being
itself connected to the problem of estimating the number of limit
cycles which can be born from polynomial perturbations of Hamiltonian
(polynomial) vector fields on the plane (see, for example, [4]).\w
{\it Remark.} The bounds are given here for a
real segment in the domain where the coefficients of the equations are
holomorphic and real on
the real axis, but they can be given in principle 
for any compact domain contained in the domain where the coefficients
are holomorphic ([3]). We also remark that the results on
which the bounds appearing here are based, are not the best
available (see discussion in [3]).\w 
Let us be given a linear differential system of equations
$$
\dot x = A(t) x
$$
with coefficients holomorphic in some simply connected complex domain.
Then the solution will be a vector valued holomorphic
function in this domain. The problem in the large is to estimate the
number of zeros the first component of the solution may have in that
domain.\w
It is known that if the coefficients of the system belong to
the class of polynomials of a given degree, a constructive
bound on the number of
zeros on some compact subset of the domain may be given which involves
just the degree and the maximum modulus of the coefficients [1].\w
In the following, another way of obtaining bounds on the number of
zeros of the solutions of such systems is discussed, which is the
consequence of a simple but somewhat surprising fact.\w
In [1], one of the aims is to bound the complex oscillation of solutions of 
polynomial
vector fields; the result for linear differential systems is obtained
as a consequence. One writes the system as a polynomial differential
system of equations with the coefficients and time treated as extra
variables, the polynomial coefficients being integral. 
Then, by subsequent differentiating, 
one obtains the derivatives of the
first variable written as polynomials with integer coefficients in
all the variables of the system. It turns out that one can effectively
bound
the moment when the ideal generated by the latter polynomials
stabilizes. This last derivative is then written as a linear
combination of the former derivatives
over the ring of polynomials  in the variables of the system, 
which is a linear scalar differential
equation for the first coordinate.
Since the coefficients of the polynomials are integers, one uses a
form of
effective Nullstellensatz to bound the coefficients of this equation
(on compact regions in the variable space). Since for linear scalar
differential equations a bound on the number of zeros exists in terms of 
the maximum of the coefficients in the region ([2]), a bound on
the number of zeros is obtained.\w
One notes that the obtained scalar differential equation has a space of
solutions of dimension much larger than the dimension of the original
system, and thus has many solutions not corresponding (at least
directly) to the solutions of the system.\w
In what we do here, we keep the integrality of the polynomial
coefficients to be able to use the effective Nullstellensatz. As in
[1], we can parameterize in this way all linear differential systems
with coefficients being polynomials of a fixed degree. In the
beginning, however, we
assume only holomorphic dependence on the time and parameters so that
not to obscure the main point.
\w
So let there be given a linear differential system of order $n$
$\dot x= A(t,p)x=(a_{ij})(t,p)(x)$, 
where $p$ stands for parameters. We assume the
dependence of the system matrix on $t,p$ to be holomorphic in some
region $U$. By $W$ we denote the projection of the domain on the
parameter space. We describe a
procedure to obtain a scalar differential  equation for the first
system variable. Differentiating and substituting the relation we have
for the first variable, we get:
$$
x_1''= \sum_{i=1}^n a_{1i}'x_i+a_{1i}x_i'=\sum_{i=1}^n a_{1i}'x_i+a_{1i}
  \sum_{i=1}^n a_{ij}x_j=\sum_{i=1}^n (a_{1i}'+\sum_{j=1}^n a_{1j}a_{ji}) x_i.
$$
Similarly we get the expressions for $x_1^{(i)}$ for any $i$. We write
$$
x_1^{(i)}=\sum a_j^{(i)}x_j= a^{(i)} x.
$$
There is a minimal $k$, 
$1\leq k\leq n$, such that $a^{(0)},a^{(1)},..,a^{(k)}$ are linearly
dependent. We then have a unique decomposition
$$
a^{(k)}=A_{k-1}a^{(k-1)}+...+A_0 a^{(0)},
$$
$A_i\in {\cal M}(U)$, 
implying that the following equation must be satisfied by 
any solution for the $x_1$ variable of the initial system
$$
x_1^{(k)}-A_{k-1}x_1^{(k-1)}-...-A_0 x_1=0.
$$
In the sequel we call this differential equation simply 'the derived
equation'. Notice that the determinant of a suitable $k$-minor of the matrix
$$
\left( a^{(0)},a^{(1)},..,a^{(k-1)}\right)
$$ 
is not
identically zero, and that multiplying by it the equation, we get an
equation with holomorphic coefficients. We thus call it the leading
coefficient of the derived equation (corresponding to the minor
chosen; of course in the case k=n there is only one such choice).\w
{\bf Claim 1.}{\it There exists an analytic set G such that for any
$p\in W-G$ any solution of the above equation is a solution
for the $x_1$ variable of the initial system.}\w
{\bf Proof.} The leading coefficient of the derived equation is a
holomorphic function in $p,t$, not identically zero. 
Thus the values of $p$ for which the leading coefficient
vanishes identically in $t$, form an analytic subset of $W$, $G$.  
For $p \in W-G$ for all but a discrete set of values of $t$ we have
$rank \ (a^{(i)}_j(t,\e)=k$.\w
Without limiting generality, we assume $U$ to be a polydisc. 
For any fixed $p_0\in W-G$, there is a natural mapping from the space of
solutions of the system to the space of solutions of the derived
equation,
which we can identify with the linear mapping 
$$
x_1(t_0,p_0),...,x_n(t_0,p_0) \ \mapsto \
x_1(t_0,p_0),...,x_1^{(k)}(t_0,p_0) 
$$
where $t_0$ is such that the leading coefficient does not vanish on
$t_0,p_0$, given by
$$
x_1^{(i)}(t_0,p_0)=(a^{(i)}_j(t_0,p_0))(x_j(t_0,p_0)).
$$
Since the leading coefficient does not vanish at $t_0,p_0$, 
$rank \ (a^{(i)}_j(t,\e))=k$, and the
mapping is surjective. Hence $\forall p_0 \in W-G$ 
every solution of the derived 
equation comes from a
solution of the system. 
\w
{\it Remark}: if the dependence is algebraic rather than holomorphic,
the set $G$ is algebraic.\n
Let us express the determinant of the $i$-th $k$-minor, the order of
the derived equation being $k$, as $b_{i0}(p)+b_{i1}(p)t+...+b_{i
k_i}t^{k_i}$. We call the ideal generated by all $b_{ij}(p)$ the
degeneracy ideal, and the corresponding variety the degeneracy
variety.\w   
We now make a small digression, and explain a basic difficulty which
might occur in
estimating the number of zeros of the first component of a solution of
the system through the derived above equation for this component. In
[2] one finds the following theorem, derived from growth estimates on
the solutions:\w
{\bf Theorem.    }{\it Let 
$$
a(t)y^{(n)}+a_{n-1}(t)y^{(n-1)}+...+a_0(t)y=0
$$
have coefficients holomorphic in a region $U$ containing a real
segment $K$. Suppose that the modulus of all $a_i(t)$ is bounded in
$U$ by $A$, and that 
$$
max_{t\in K}|a(t)| \ \geq \ a.  
$$
Then the number of zeros on $K$ of any holomorphic in $U$ solution of the
equation, real on the real axis, is not greater than 
$$
(A/a+n)^\mu,
$$
where $\mu$ depends only on the geometry of $U$ and $K$.}\w
The main obstacle in applying this theorem to our case lies in the
fact that, as parameters vary, $max_{t\in K}|a(t,p)|$ (which depends
on the parameters) {\it a priori} 
may become arbitrarily small - the situation of
singular perturbations - and the estimates explode.\w 
Our next goal is to show that the scalar differential equation we
derived for $x_1$
cannot be in fact singularly perturbed, at least when one parameter is
involved. This is somewhat unexpected,
and we explain why. Consider the following very simple
$\e$-dependent system:
$$
\dot x = x+\e y \ \ \ \ \dot y= x+y
$$
One might expect that when $\e=0$, the derived equation, which is of
second order, will simply
degenerate into the first order $\dot x=x$, implying that the derived
equation is singularly perturbed. However, $\ddot x=x+\e
y+\e(x+y)=(1+\e)x+2\e y$, and the derived equation is
$$
\ddot x -2\dot x+(1-\e)x=0.
$$
This equation is not singularly perturbed. It has the property that
for the exceptional value of the parameter, an additional solution
appears which does not come from the original system, namely $te^t$.
To explain why the singular perturbations cannot occur also in the
general case, we proceed as follows. Until further notice the dimension of the
parameter space is one, and the parameter will be denoted by $\e$.\w
{\bf Definition.} We call a meromorphic function $f(t,\e)$ singularly
perturbed at $\e_0$, if when expressed as a quotient of two holomorphic
functions $f(t,\e)=u(t,\e)/v(t,\e)$, $u(t,\e)$ is divisible by a higher or
equal power of $(\e-\e_0)$ than the one that divides $v(t,\e)$.\w
{\bf Definition.} A (scalar) differential equation with meromorphic
coefficients depending on one parameter, normalized so that the
leading coefficient is equal to 1, will be called singularly perturbed if
at least one of its coefficients is singularly perturbed at some point
in the parameter space.\w
We make the following obvious claim.\w
{\bf Lemma 1. }{\it Let us be given a differential
equation with $t,\e$-dependent coefficients meromorphic in a domain $U$, 
singularly perturbed at $\e_0$.
Then there is a point $t_0$, such that in the neighborhood of
$(t_0,\e_0)\in U$, the equation takes the form
$$
(\e-\e_0)^s y^{(n)}+a_{n-1}(t,\e) y^{(n-1)}+...+a_0(t,\e)y =0,
$$
$s>0$, with coefficients holomorphic in the neighborhood of
$(t_0,\e_0)$, and at least one
of them is not identically zero when $\e=\e_0$.}\w
{\bf Theorem 1.} {\it Let a differential equation of order $n$ with
coefficients which depend meromorphically on $t,\e$ in some polydisc
$U\times V$, be 
singularly perturbed for at least one parameter value $\e_0$, 
and let it have l solutions $y_1(t,\e),...,y_l(t,\e)$
holomorphic in $t$ and $\e$ in that polydisc.   
Then for each fixed parameter value $\hat \e$ from $V$,}
$$
dim_{\C} \ span_{\C}\left( y_1(t,\hat \e),...,y_l(t,\hat \e)\right) \ \leq \
\ n-1.
$$
{\bf Proof.} Using Lemma 1, we write the equation in the form given
there in some smaller polydisc, and assume $t_0=0, \ \e_0=0$:
$$
\e^s y^{(n)}+a_{n-1}(t,\e) y^{(n-1)}+...+a_0(t,\e)y =0,
$$  
with $a_k(0,0)\neq 0$. We claim that if two holomorphic in $t,\e$
solutions are such that they have equal
$y^{(n-1)}(0,\e),...,y^{(k+1)}(0,\e),y^{(k-1)}(0,\e),...,y(0,\e)$,
they are in fact identical. Indeed, otherwise we get that there exists
a nonzero holomorphic in $t,\e$ solution $y(t,\e)$ such that
$$
y^{(n-1)}(0,\e)=...=y^{(k+1)}(0,\e)=y^{(k-1)}(0,\e)=...=y(0,\e)=0.
$$
Since $y(t,\e)$ is not identically zero, $y^{(k)}(0,\e)$ cannot
be identically zero. Thus  $y^{(k)}(0,\e)=\e^h
g(\e)$, $h\geq 0, \ g(0)\neq 0$. Assume $h>0$. 
$y(t,0)$ is the solution of the differential equation
with $\e$ set to zero. Since $a_k(0,0)\neq 0$ and the initial
conditions at $\e=0$ are zero ($h>0$), we have $y(t,0)=0$
identically. This can only happen if $y(t,\e)=\e y_1 (t,\e)$,
$y_1(t,\e)$ being also a holomorphic solution of the system with the
same zero initial conditions except for $y_1^{(k)}(0,\e)=\e^{h-1}g(\e)$. We
continue in the same way, until we get a holomorphic solution
$y_h(t,\e)=y(t,\e)/\e^h$ with the initial conditions
$$
y^{(n-1)}(0,\e)=...=y^{(k+1)}(0,\e)=y^{(k-1)}(0,\e)=...=y(0,\e)=0,
$$
and $y^{(k)}(o,\e)=g(\e)$, $g(0)\neq 0$. However substituting
$y_h(t,\e)$ into the original equation and setting $t=0$, we get:
$$
\e^s y^{(n)}(0,\e)+a_{k}(0,\e)y^{(k)}(0,\e)=0, 
$$
which is a contradiction ($a_k(0,0),y^{(k)}(0,0)\neq 0$).\w
Thus in the smaller polydisc the holomorphic in $t,\e$ solutions
can be identified with sets of holomorphically dependent on $\e$ sets
of initial conditions
$y^{(n-1)}(0,\e),...,y^{(k+1)}(0,\e),y^{(k-1)}(0,\e),...,y(0,\e)$,
omitting $y^{(k)}(0,\e)$. 
Now we prove the theorem. Suppose we have holomorphic in $U\times V$
solutions of our singularly perturbed equation 
$y_1(t,\e),..,y_l(t,\e)$. Over the field meromorphic in $V$ functions
there can be at most $n-1$ linearly independent such solutions since
the latter are identified with ${\cal O}(V)^{n-1}$. So we choose a
maximal set $S$ of linearly independent over ${\cal M}(V)$ solutions and
express all the other as their linear combination over ${\cal M}(V)$. It
follows that all the other solutions lie in the space spanned by the
solutions from $S$ for all but a discrete set of $\e$-s. Thus for
almost all values of $\e$, the space of holomorphic in $t,\e$ solutions
is $n-1$ dimensional, which implies that it is at most $n-1$
dimensional for all values of $\e$ by the semicontinuity of the dimension.\w
{\bf Theorem 2. }{\it Suppose that we are given a linear differential
system with coefficients depending on one
parameter, holomorphic both in $t$ and $\e$ in a polydisc $U\times V$. 
Then the derived scalar differential equation for the
variable $x_1$ cannot be singularly perturbed.}\w
{\bf Proof.} Returning to notations of proof of Claim 1, we take a
point $(t_0,\e_0)\in U-Z$, in which the leading coefficient of the
derived equation is not zero, and $rank \ (a^{(i)}_j(t,\e)=k$.\w
For each $m$, there is a holomorphic in $U$ solution of the system
defined by the initial conditions (constant in $\e$):
$$
x_1(t_0,\e)=0 \ ... \ x_m(t_0,\e)=1 \ ... \ x_n(t_0,\e)=0.
$$
We denote the 1-st component of this solution by $x_1^m(t,\e)$. By
discussion from the proof of Claim 1, since
$x_1^1(t,\e_0),...,x_1^n(t,\e_0)$ span the space of solutions of the
derived equation at $\e_0$, and since at this (generic) point the
order of the derived equation is $k$, we have
$$
dim_\C \ span_{\C} \ x_1^1(t,\e_0),...,x_1^n(t,\e_0) \ = \ k.
$$
However, if the derived equation is singularly perturbed even at only
one point $\e_1\in V$, from Theorem 1 it would follow that for every
$\e_0\in V$
$$
dim_\C \ span_{\C} \ x_1^1(t,\e_0),...,x_1^n(t,\e_0) \ \leq \ k-1,
$$
giving a contradiction. Thus the derived equation cannot be singularly
perturbed.\w
We explain now our goal in what follows. As it was already mentioned,
we can parameterize the problem of general linear system of order $n$ with
polynomial coefficients of degree $d$ as:
$$
\dot x_i=\sum_{j=1}^n \sum_{k=0}^d a_{ijk} t^k x_j \ \ \ i=1,..,n
$$
the parameters being $a_{ijk}$, thus making the system matrix to be 
polynomial in the time and parameters, with integer coefficients.
Ideally we would like to bound the number of zeros which the first
variable of this system might have, but the results available here are
still not permitting such general bounds. Instead, we restrict
ourselves to the case when the system matrix is polynomial in time and
{\it one} additional parameter, then making several remarks about the
general linear differential system with polynomial coefficients as
above.\n
So let us have a
linear differential system of order $n$, with coefficients being
polynomials in $t,\e$ of degree $d$, with integer coefficients of
maximum modulus $M$. Let the order of the derived equation be $k$.\w
{\bf Lemma 2. } {\it Let a nonsingularly perturbed rational function
be written as
$$
{a_0(\e)+a_1(\e)t+...+a_{k}(\e)t^{k} \over 
b_0(\e)+b_1(\e)t+...+b_{l}(\e)t^{l}}.
$$
Then each $a_j$ belongs to the ideal generated by
$b_0,b_1,..,b_{l}$.}\w
{\bf Proof.} $a_j$ must vanish at points where all
$b_1,..,b_{l}$ vanish (otherwise the fraction would be singularly
perturbed). Write $b_i(\e)=c_i(\e)b(\e)$, where $b(\e)$ is the common
factor. Then $c_i(\e)$ have no common zero and therefore, since we are
working over $\C$, $1=h_1(\e)c_1(\e)+..+h_l(\e)c_l(\e)$, for some
polynomials $h_1,..,h_l$. But $a_j$ must be divisible by $b(\e)$:
otherwise we get singular perturbations. Thus:
$$
a_j(\e)=d(\e)b(\e)=d(\e)b(\e)(h_1(\e)c_1(\e)+..+h_l(\e)c_l(\e))=
d(\e)h_1(\e)b_1(\e)+..+d(\e)h_l(\e)b_l(\e).
$$
We derive now some computational estimates we will need.\w
{\bf Lemma 3.} {\it Suppose we are given a non-singularly perturbed
rational function of $t,\e$, 
$$
{a_0(\e)+a_1(\e)t+...+a_{k}(\e)t^{k} \over 
b_0(\e)+b_1(\e)t+...+b_{l}(\e)t^{l}},
$$
with degrees of the numerator and the denominator being at most $d$,
and the coefficients being integers of maximum modulus M (we mean by
this possibly complex integers). Then one may write:
$$
{a_0(\e)+a_1(\e)t+...+a_{k}(\e)t^{k} \over 
b_0(\e)+b_1(\e)t+...+b_{l}(\e)t^{l}} = {\sum_{i=1}^k t^i \sum_{j=1}^l
h_{ij}(\e) b_j(\e) \over b_0(\e)+b_1(\e)t+...+b_{l}(\e)t^{l}},
$$
with $h_{ij}(\e)$ being polynomials of degree at most $2d-1$ with
integer coefficients of maximum modulus $(eM)^{Cd^3}$ for some
constant $C$.}\w
{\bf Proof.} From Lemma 2, $a_i=\sum h_{ij}b_j$. 
By [5], since $\forall i \ deg(a_i)\leq d$, $\forall j \
deg(b_j)\leq d$, we have that $a_i=\sum h_{ij}b_j, \ deg(h_{ij})\leq
2d-1$. Observe that this gives us an integer coefficient linear system
of dimension at most $2d(d+1)$, maximum modulus of the coefficients not
greater than $M$. By the integrality trick this gives us the estimate
$(2d(d+1))!M^{2d(d+1)}\leq (eM)^{Cd^3}$ 
on the coefficients of the polynomials $h_{ij}$.\w
{\bf Lemma 4.} {\it In the setting of Lemma 3, suppose that $|\e|$ is
bounded by $E$. Then the above fraction can always be written as
$$
{\sum_{i=1}^d   c_i t^i \over \sum_{i=1}^d d_i t^i}
$$
with $c_i\leq (eM)^{Cd^3}(E^{2d-1}+1)$ for some constant $C$
, with some 
$d_i$ equal to $1$ and all $|d_i|$ being not greater than 1.}\w
{\bf Proof.} For each $\e$ in the disc, take the $b_i$ having the
maximum modulus, and divide by it both the numerator and the
denominator: use Lemma 3 (the constant here is different from the one
in Lemma 3).\w  
{\bf Lemma 5.} {\it Suppose we are given a non-singularly perturbed
differential equation of order $k$
$$
y^{(k)}+{a_{k-1}(t,\e)\over b(t,\e)} y^{(k-1)}+...+{a_0(t,\e)\over
b(t,\e)}y =0 ,
$$ 
all $a_i(t,\e)$ and $b(t,\e)$ being polynomials
of degree at most $d$, having integer coefficients of
maximum modulus $M$.\n
Let $[-R/2,R/2]$ be a segment on the real axis containing $[-1,1]$. 
Fix any real $\e$ of modulus less
than $E$. Let $f$ 
be a holomorphic solution of the equation for this parameter value, 
real on the real axis. Then
the number of zeros $f$ can have on $[-R/2,R/2]$, is bounded by
$$
\left( (Me)^{Cd^3}(E^{2d}+1)R^{d+1}+k\right)^{\sigma},
$$
with $\sigma$ and C being some constants.}\w
{\bf Proof.}
By Lemma 4, we can write the equation in the form
$$
(t^s+\sum_{i\neq s}d_i t^i)y^{(k)}+(c_{ k-1 \ 0}+..+c_{
k-1 \ d}t^d)y^{(k-1)}+...+ (c_{0 0}+...+c_{0 d}t^d)y=0.
$$
with $|c_{ij}|$ being bounded by the bound from Lemma 2, and all
$|d_i|\leq 1$. If we could
find now some point in $[-R/2,R/2]$ where $t^s+\sum_{i\neq s}d_i t^i$ has some
value bounded away from zero, we are done. This because
then we may just apply the Theorem cited from [2].\n
So we write
$$
|t^s+\sum_{i\neq s}d_i t^i|\geq
|t^s+d_{s-1}t^{s-1}+...+d_0|-|t^{s+1}(d_{s+1}+d_{s+2}t+..)|
$$
and notice that by Cartan lemma ([2]), one may delete disks of total
diameter $h$ from $\C$, such that on the rest of $\C$ the modulus of
the monic polynomial $t^s+d_{s-1}t^{s-1}+...+d_0$ is $\geq
(h/4e)^s$. Thus for each $h>0$ there exists a point $t_h$ of 
modulus at most $h/2$, at which
the value of the polynomial is $\geq (h/4e)^s$. But then:
$$
|t_h^s+\sum_{i\neq s}d_i t_h^i|\geq
|t_h^s+d_{s-1}t_h^{s-1}+...+d_0|-|t_h^{s+1}(d_{s+1}+d_{s+2}t_h+..)|\geq
(h/4e)^s-h^{s+1}(d-s). 
$$
Taking for example $h=1/(2(4e)^s(d-s))$ we get that there exists a point
on the unit segment at which $|t^s+\sum_{i\neq s}d_i t^i|$ is bounded
from below by 
$$
{1\over (4e)^{s^2+s} 2^{s+1} (d-s)^s}\geq {1\over (4e)^{s^2+s} 2^{s+1} d^s}
$$
As $s$ can be from $0$ to $d-1$ (it can be shown that the case $s=d$ is
also covered by the obtained estimate), we are led to the conclusion
that there will be a point in $[-R/2,R/2]$ where 
$$
|t^s+\sum_{i\neq s}d_i t^i|\geq {1\over (4e)^{d^2+d} 2^{d+1} d^d}.
$$
In turn, $|\sum_{j=0}^d c_{ij}t^j|$ cannot, when $\e$ varies over disk
of radius $E$ and $t$ over $D_R$, be larger than
$$
\sum_{j=0}^d R^j\cdot (eM)^{Cd^3}(E^{2d-1}+1) \ \ \ \leq \ \ \ 
(eM)^{Cd^3} (E^{2d-1}+1)R^{d+1}.
$$
By Theorem cited from [2], one deduces that for any given real $|\e|\leq E$,
any holomorphic in $t$ 
solution of the given differential equation, real on the real
axis, may have not more than
$$
\left({ (eM)^{Cd^3} (E^{2d-1}+1)R^{d+1} \over {1/
(4e)^{d^2+d} 2^{d+1} d^d}}+k\right)^{\sigma}
$$
zeros on $[-R/2,R/2]$, $\sigma$ not depending on $R$ because of properties of
$\sigma(K,U)$ from the Theorem. Simplifying, one gets the bound we stated.\w
To estimate the number of zeros a solution (more precisely, its $x_1$
component) of the initial linear system
$$
\dot x=A(t,\e)x
$$
might have, we now have only to bound the maximum modulus of
the (integer) coefficients of the derived equation. We will do the
computation for an arbitrary number of parameters; to emphasize this
we shall write $p$ instead of $\e$. To remind the notations, we denote the
dimension of the system by $n$, the degree of the polynomial
coefficients by $d$, the maximum modulus of the coefficients of those
polynomials by $M$, the number of parameters by $q$, 
and the order of the derived equation for $x_1$ by
$k$.\w
{\bf Lemma 6.} {\it The following recurrent formulas hold for $a^{(i)}$,
defined in the beginning of this text:}
$$
a^{(i)}=\dot a^{(i-1)} +A(t,p)^T a^{(i-1)}.
$$
{\bf Proof.} Immediate.\w
The following is a computation:\w
{\bf Lemma 7.} {\it Let $P_1$, $P_2$ be two polynomials of degrees
$d_1,d_2$ in s variables, with the modulus of their coefficients being
not larger than $M_1,M_2$, respectively. 
Then the degree of $P_1 P_2$ is $d_1+d_2$ and the coefficients are not
larger in absolute value than $(1+min(d_1,d_2))^sM_1M_2$.}\w
{\bf Lemma 8.} {\it The degrees and the maximum modulus of
coefficients of the entries of $a^{(i)}$ are $di$ and 
$n^i (d+(d+1)^{q+1}M)^i$, respectively.}\w
{\bf Proof.} By writing recursive relations.\w
We now compute the maximum modulus and the degrees for
the coefficients of the derived equation.\w
{\bf Lemma 9.} Let 
$$
y^{(k)}+{\a_{k-1}(t,p)\over \b(t,p)} y^{(k-1)}+...+{\a_0(t,p)\over
\b(t,p)}y =0 
$$ 
be the equation derived for $x_1$ from the original system. Then the
degrees of both $\a_i$ and $\b$ are bounded by $k(k+1)d/2$, and the maximum
modulus for the (integer) coefficients of those polynomials is not
larger than $(eM(d+1))^{C(q+1)n^3}$ for some constant $C$.\w
{\bf Proof.} Again by recursive relations.\w
{\bf Theorem 2.} {\it Let the initial linear differential system have as
coefficients polynomials in $t,\e$ ($\e$ one dimensional) with integer
coefficients of maximum modulus $M$, and let $[-R/2,R/2]$ be 
a segment on the real axis containing $[-1,1]$. Fix any real 
$\e$ of modulus less than $E$. Let $f$ 
be a holomorphic solution of the equation for this parameter value, 
real on the real axis. Then
the number of zeros $f$ can have on $[-R/2,R/2]$, is bounded by
$$
\left( (Me)^{Cn^9 d^4}(E^{n(n+1)d}+1)R^{n(n+1)d/2+1}+n\right)^{\sigma},
$$
for some constants $C$ and $\sigma$.}\w
{\bf Proof.} Just substitute the bounds stated in Lemma 9 with $q=1$ 
into Lemma 5 and simplify the expression.\w
Of course, the concrete form of the bound given above is irrelevant; 
what is important is the rate of its growth, which
is simple exponential in $n,d$ compared with the tower of exponents 
which can be derived from [1] in this situation.\n 
We now make some remarks on the general case when many parameters are 
allowed. First, it is clear that
the absence of singular perturbations in the case of dependence on one
parameter implies absence of singular perturbations 
along any given holomorphic curve in the parameter
space, and the following is true.\w
{\bf Lemma 2a. } {\it Let a rational function dependent on several
parameters be written as
$$
{a_0(p)+a_1(p)t+...+a_{k}(p)t^{k} \over 
b_0(p)+b_1(p)t+...+b_{l}(p)t^{l}},
$$
and let it be such that it is not singularly perturbed along any line
in the parameter space.
Then each $a_j$ vanishes on the variety where all
$b_0,b_1,..,b_{l}$ vanish.}\w 
We cannot in general have the conclusion of Lemma 2  
as the example of the nonsingularly perturbed fraction
$$
ab\over a^2+b^2 t
$$
shows. That is, in the derived equation for system with polynomial
coefficients depending on many parameters we can only guarantee that
the coefficients in the numerator vanish on the degeneracy variety,
but not that they in fact belong to the degeneracy ideal. However, if
by some additional reasoning
for, say, the general linear system with polynomial coefficients
$$
\dot x_i=\sum_{j=1}^n \sum_{k=0}^d a_{ijk} t^k x_j \ \ \ i=1,..,n,
$$ 
it will turn out to be true, then the
bound, though drastically grows because of the double exponential
in $n^2(d+1)$ bounds in the effective Nullstellensatz (see [5]), has
the asymptotic form (in terms of Theorem 2)
$$
e^{e^{e^{P(n,d)}}}(ER)^{Cn^2 d},
$$
where $P$ is a polynomial. This is still a considerable
improvement compared to [1] (of course if the assumption is valid).\w
\begin{center}{\bf \ \ \ \ \ \ \ References.}\end{center}
1.  Novikov, D.; Yakovenko, S. {\it Trajectories of polynomial vector fields and
ascending chains of polynomial ideals.} Ann. Inst. Fourier (Grenoble) 49 (1999), no. 2, 563--609.\w  
2. Il'yashenko, Y.; Yakovenko, S. {\it Counting real zeros of
analytic functions satisfying linear ordinary differential equations.} J. Differential Equations 126
(1996), no. 1, 87--105.\w
3. Yakovenko S. {\it On functions and curves defined by ordinary differential equations}
, to appear in Proceedings of the Arnoldfest (Ed. by Bierstone, Khesin, Khovanskii, Marsden), Fields Institute Communications, 1998.\w
4. Il'yashenko, Y.; Yakovenko, S. {\it Double exponential
estimate for the number of zeros of complete abelian integrals and rational envelopes of linear
ordinary differential equations with an irreducible monodromy group.} Invent. Math. 121 (1995), no. 3,
613--650.\w
5. Shiffman, B.  
{\it Degree bounds for the division problem in polynomial ideals.} 
Michigan Math. J. 36 (1989), no. 2, 163--171. 

\end{document}